\documentclass[final]{dmtcs-episciences}

\usepackage[utf8]{inputenc}
\usepackage{amsmath,amssymb,mathrsfs,amsthm}

\newtheorem{theorem}{Theorem}[section]
\newtheorem{lemma}{Lemma}[section]
\newtheorem{proposition}{Proposition}[section]
\newtheorem{corollary}[theorem]{Corollary}
\newtheorem{fact}{Fact}
\newtheorem{example}{Example}[section]

\title[Edge-disjoint Steiner trees in Cartesian product networks]
{Constructing edge-disjoint Steiner trees in Cartesian product networks}

\author[Li et al.]{%
Rui Li\affiliationmark{1}
\and Gregory Gutin\affiliationmark{2}
\and He Zhang\affiliationmark{3}\\
\and Zhao Wang\affiliationmark{4}
\and Xiaoyan Zhang\affiliationmark{5}
\and Yaping Mao\affiliationmark{6}\thanks{Corresponding author.}}

\affiliation{%
School of Computer Science and Technology, Xidian University, Xi'an, China\\
Department of Computer Science, Royal Holloway, University of London, London, UK\\
School of Mathematics and Statistics, Qinghai Normal University, Xining, China\\
College of Science, China Jiliang University, Hangzhou, China\\
School of Mathematical Science and Institute of Mathematics, Nanjing Normal University, Nanjing, China; Ministry of Education Key Laboratory of NSLSCS\\
Ministry of Education Key Laboratory of NSLSCS, School of Mathematics and Statistics, and Academy of Plateau Science and Sustainability, Qinghai Normal University, Xining, China}


\hypersetup{
  pdftitle={Constructing edge-disjoint Steiner trees in Cartesian product networks},
  pdfauthor={Rui Li, Gregory Gutin, He Zhang, Zhao Wang, Xiaoyan Zhang, Yaping Mao},
  pdfkeywords={tree, Steiner tree, Cartesian product, generalized connectivity}
}

\begin{document}
\publicationdata{vol. 28:2}{2026}{15}{10.46298/dmtcs.10921}{2023-02-07; 2023-02-07; 2024-05-11}{2026-02-12}
\maketitle

\begin{abstract}
Cartesian product networks are often used to combine two given networks with
established properties into a new network that inherits properties from both.
For a graph $F=(V,E)$ and a set $S\subseteq V(F)$ of at least two vertices, an
$S$-Steiner tree (or simply an $S$-tree) is a tree subgraph $T=(V',E')$ of $F$
with $S\subseteq V'$.  The generalized local edge-connectivity $\lambda(S)$ is
the maximum number of pairwise edge-disjoint $S$-Steiner trees in $F$.  For an
integer $k$ with $2\leq k\leq |V(F)|$, the generalized $k$-edge-connectivity of
$F$ is
\[
  \lambda_k(F)=\min\{\lambda(S): S\subseteq V(F),\ |S|=k\}.
\]
We establish sharp upper and lower bounds for $\lambda_k(G\Box H)$, where
$\Box$ denotes the Cartesian product and $G$ and $H$ are graphs.
\end{abstract}

\section{Introduction}

A processor network can be represented by a graph, in which vertices
represent processors and edges represent communication links. The process of
sending a message from a source node to all other nodes in a network is called
broadcasting. This can be accomplished by having each node repeatedly receive
and forward messages. Since some nodes and links may be faulty, multiple
copies of a message can be disseminated through edge-disjoint Steiner trees.
A broadcast succeeds if the healthy nodes obtain the correct message from the
source node within a given time. Considerable attention has been devoted to
fault-tolerant broadcasting in networks \cite{Fragopoulou, Hedetniemi, Jalote, Ramanathan}. To measure fault tolerance, path structures connecting two nodes must be
generalized to tree structures connecting more than two nodes; see \cite{KuWangHung, QinHaoChang,
ZhaoHaoWu, ZhaoHaoCheng, LiTuYu, ZhaoHao, LiMaYangWang, LinZhang,
LiWang}.

All graphs considered in this paper are undirected, finite and
simple. We refer to the book \cite{BondyMurty} for graph-theoretic
notation and terminology not defined here. For a graph $F$, let
$V(F)$, $E(F)$ and $\delta(F)$ denote the set of vertices, the set
of edges and the minimum degree of $F$, respectively.
Edge-connectivity is one of the most basic concepts of
graph theory, both in the combinatorial and algorithmic senses. It is well-known that the classical
edge-connectivity has two equivalent definitions. The \emph{edge-connectivity} of $F$, written $\lambda(F)$, is the minimum size
of an edge-subset $M\subseteq E(F)$ such that $F\setminus M$ is
disconnected. We call this the ``cut'' definition of edge-connectivity. The ``path''
definition is as follows. First, define
$\lambda_F(x,y)$ for a pair of distinct vertices $x$ and $y$ of $F$,
which is the maximum number of edge-disjoint paths between $x$ and
$y$ in $F$. This parameter is called the \emph{local
edge-connectivity} of $x$ and $y$ in $F$. For a graph $F$, define the global quantity $\lambda^*(F)=\min\{\lambda_F(x,y)\,|\, x, y\in
V(F), \, x\neq y\}$. Menger's theorem says that $\lambda(F)$ is
equal to $\lambda^*(F)$. This result can be found in any standard textbook; see, for example,
\cite{BondyMurty}.

Although there are many elegant results on edge-connectivity, the classical
notion may not be general enough to capture some computational settings. This
motivates its generalization. The generalized connectivity of a graph,
introduced by Hager \cite{Hager}, is a natural generalization of the ``path''
definition of connectivity. For a graph
$F=(V,E)$ and a set $S\subseteq V(F)$ of at least two vertices,
\emph{an $S$-Steiner tree} or \emph{a Steiner tree connecting $S$}
(or simply, \emph{an $S$-tree}) is a subgraph $T=(V',E')$ of $F$
that is a tree with $S\subseteq V'$. When $|S|=2$, a minimal
$S$-Steiner tree is just a path connecting the two vertices of $S$.
Li et al.\ introduced the corresponding concept of generalized
edge-connectivity in \cite{LMS}. For $S\subseteq V(F)$ with $|S|\geq 2$, the \emph{generalized local edge-connectivity} $\lambda(S)$ is the maximum
number of edge-disjoint Steiner trees connecting $S$ in $F$. For an
integer $k$ with $2\leq k\leq |V(F)|$, the \emph{generalized
$k$-edge-connectivity} $\lambda_k(F)$ of a graph $F$ is defined as
$\lambda_k(F)=\min\{\lambda(S):S\subseteq V(F),\ |S|=k\}$.
When $|S|=2$, $\lambda_2(F)$ is just the standard edge-connectivity
$\lambda(F)$ of $F$, that is, $\lambda_2(F)=\lambda(F)$, which explains the term generalized edge-connectivity. Also set $\lambda_k(F)=0$ when $F$ is
disconnected. Table~\ref{tab:connectivity-comparison} summarizes this generalization.
\begin{table}[!h]
\centering
\small
\resizebox{\linewidth}{!}{%
\begin{tabular}{|c|c|}
\hline
Edge-connectivity & Generalized edge-connectivity\\
\hline
$S=\{x,y\}\subseteq V(F)$, $|S|=2$
  & $S\subseteq V(F)$, $|S|\geq 2$\\
\hline
$\left\{
\begin{array}{l}
\mathcal{P}_{x,y}=\{P_1,P_2,\ldots,P_{\ell}\},\\
\{x,y\}\subseteq V(P_i),\\
E(P_i)\cap E(P_j)=\varnothing\quad(i\ne j)
\end{array}
\right.$
&
$\left\{
\begin{array}{l}
\mathcal{T}_{S}=\{T_1,T_2,\ldots,T_{\ell}\},\\
S\subseteq V(T_i),\\
E(T_i)\cap E(T_j)=\varnothing\quad(i\ne j)
\end{array}
\right.$\\
\hline
$\lambda_F(x,y)=\max |\mathcal{P}_{x,y}|$
  & $\lambda(S)=\max |\mathcal{T}_{S}|$\\
\hline
$\lambda(F)=\displaystyle\min_{\substack{x,y\in V(F)\\x\ne y}}\lambda_F(x,y)$
  & $\lambda_k(F)=\displaystyle\min_{\substack{S\subseteq V(F)\\|S|=k}}\lambda(S)$\\
\hline
\end{tabular}%
}
\caption{Classical edge-connectivity and generalized edge-connectivity.}
\label{tab:connectivity-comparison}
\end{table}

Zhao and Hao \cite{ZhaoHao} studied the generalized connectivity of
alternating group graphs and $(n,k)$-star graphs. Li et al.~\cite{LiTuYu} determined the exact values of generalized $3$-connectivity
of star graphs and bubble-sort graphs. Zhao et al.~\cite{ZhaoHaoWu}
investigated the generalized $3$-connectivity of some regular
networks, and Zhao et al.~\cite{ZhaoHaoCheng} studied the
generalized connectivity of dual cubes. Li and Wang \cite{LiWang}
obtained the exact values of generalized $3$-(edge)-connectivity of
recursive circulants. There are many results on generalized
edge-connectivity; see the book \cite{LMbook} and papers
\cite{LiWang, LiLiSun, LiMaYangWang, LiMaoExtremal, LiMaoProduct,
LiMaoNordhaus, LiTuYu, LinZhang, Sun}.

Being a natural combinatorial measure, generalized $k$-connectivity
can be motivated by its interesting interpretation in practice. For
example, suppose that $F$ represents a network. To connect a pair of vertices of $F$, one uses a path. To connect a set $S$
of at least three vertices, one instead uses a tree.
Such a tree is called a \emph{Steiner tree} and is widely used in the physical design of VLSI circuits (see
\cite{Grotschel, GrotschelMartinWeismantel, Sherwani}). In this
application, a Steiner tree is needed to share an electric signal by
a set of terminal nodes. Steiner trees are also used in computer
communication networks (see \cite{Du}) and optical wireless
communication networks (see \cite{Cheng}). The number of mutually independent ways to connect the terminals provides a
measure of the robustness of the network. The generalized $k$-connectivity can be used to measure the capability of a network $F$ to connect any $k$
vertices in $F$.

A set of spanning trees in a graph $F$ is said to be
\emph{independent} if all the trees are rooted at the same vertex $r$
such that, for any other vertex $v\neq r$ in $F$, the paths from
$v$ to $r$ in any two trees have no common vertex except the two
endvertices $v$ and $r$. Itai and Rodeh \cite{ItaiRodeh} first introduced
the concept of ISTs to investigate the reliability of distributed
networks. For more details, we refer to \cite{ItaiRodeh,
YangChangPaiChan, YangChangTangWang}.

Product networks are based upon the idea of using the product as a
tool for combining two networks with established properties to
obtain a new one that inherits properties from both
\cite{DayAl-Ayyoub}. There has been an increasing interest in a
class of interconnection networks called Cartesian product networks;
see \cite{DayAl-Ayyoub, Hammack, KuWangHung}.

The \emph{Cartesian product} of two graphs $G$ and $H$, written as
$G\Box H$, is the graph with vertex set $V(G)\times V(H)$, in which
two vertices $(u,v)$ and $(u',v')$ are adjacent if and only if
$u=u'$ and $vv'\in E(H)$, or $v=v'$ and $uu'\in E(G)$.

Xu and Yang obtained a formula for the edge-connectivity of a
Cartesian product in 2006.
\begin{theorem}[\cite{XuYang}]\label{Thm:Ramseyn-24m}
Let $G$ and $H$ be graphs on at least two vertices. Then
\[
\lambda(G\Box H)=\min \{\lambda(G) |V(H)|,\lambda(H)
|V(G)|,\delta(G)+\delta(H)\}.
\]
\end{theorem}

For $\lambda_3(G\Box H)$, Sun in \cite{Sun} obtained the following
lower bound.
\begin{theorem}[\cite{Sun}] \label{th7-1-10}
Let $G$ and $H$ be connected graphs. Then
\[
\lambda_3(G\Box H)\geq \lambda_3(G)+\lambda_3(H).
\]
Moreover, the lower bound is sharp.
\end{theorem}

Ku, Wang, and Hung \cite{KuWangHung} obtained the following result.

\begin{theorem}[\cite{KuWangHung}]\label{th7-2-8}
Let $G$ be a graph having $\ell_1$ edge-disjoint spanning trees and
$H$ be a graph having $\ell_2$ edge-disjoint spanning trees. Then
the product network $G\Box H$ has $\ell_1+\ell_2-1$ edge-disjoint
spanning trees.
\end{theorem}

The following corollary is immediate.
\begin{corollary}[\cite{KuWangHung}]\label{cor:ku-wang-hung}
Let $G$ and $H$ be connected graphs of orders $n_G$ and $n_H$,
respectively. Then
\[
\lambda_{n_G n_H}(G\Box H)\geq
\lambda_{n_G}(G)+\lambda_{n_H}(H)-1.
\]
Moreover, the lower bound is sharp.
\end{corollary}

In this paper, we derive the following lower bound.
\begin{theorem}\label{th-Lower}
Let $k$ be an integer with $k\geq 2$, and let $G$ and $H$ be two
connected graphs with at least $k+1$ vertices. If $\lambda_k(H)\geq
\lambda_k(G)\geq \lfloor \frac{k}{2}\rfloor$, then
\[
\lambda_k(G\Box H)\geq \lambda_k(G)+\lambda_k(H).
\]
Moreover, the lower bound is sharp.
\end{theorem}

The following upper and lower bounds, due to Li and Mao, are in
\cite{LiMaoNordhaus}.
\begin{proposition}[\cite{LiMaoNordhaus}]\label{pro1-1}
For a connected graph $G$ of order $n$ and $3\leq k\leq n$,
\[
\frac{1}{2}\lambda(G)\leq \lambda_k(G)\leq \lambda(G).
\]
Moreover, the lower bound is sharp.
\end{proposition}

\begin{proposition}\label{pro-Upper}
Let $k$ be an integer with $k\geq 2$, and let $G$ and $H$ be two
connected graphs with at least $k+1$ vertices. Then
\[
\lambda_k(G\Box H)\leq
\min\{2\lambda_k(G)|V(H)|,2\lambda_k(H)|V(G)|,\delta(G)+\delta(H)\}.
\]
Moreover, the upper bound is sharp.
\end{proposition}
\begin{proof}
From Proposition~\ref{pro1-1}, we have $\lambda(G)\leq
2\lambda_k(G)$, and hence
\begin{align*}
\lambda_k(G\Box H)\leq \lambda(G\Box H) &\leq \min \{\lambda(G)
|V(H)|,\lambda(H)
|V(G)|,\delta(G)+\delta(H)\}\\[0.2cm]
&\leq \min \{2\lambda_k(G) |V(H)|,2\lambda_k(H)
|V(G)|,\delta(G)+\delta(H)\},
\end{align*}
as required.
\end{proof}

\section{Lower bounds}

Throughout this section, let $G$ and $H$ be two connected graphs with
$V(G)=\{u_1,u_2,\ldots,u_{n}\}$ and $V(H)=\{v_1,v_2,\ldots,v_{m}\}$,
respectively. Then $V(G\Box H)=\{(u_i,v_j)\,|\,1\leq i\leq n, \
1\leq j\leq m\}$. For $v\in V(H)$, let $G(v)$ denote the
subgraph of $G\Box H$ induced by the vertex set $\{(u_i,v)\,|\,1\leq
i\leq n\}$. Similarly, for $u\in V(G)$, let $H(u)$ denote the
subgraph of $G\Box H$ induced by the vertex set $\{(u,v_j)\,|\,1\leq
j\leq m\}$.

Let
$S=\{(u_{i_1},v_{j_1}),(u_{i_2},v_{j_2}),\ldots,(u_{i_k},v_{j_k})\}\subseteq
V(G\Box H)$, $S_G=\{u_{i_1},u_{i_2},$ $\ldots,u_{i_k}\}$, and
$S_H=\{v_{j_1},v_{j_2},\ldots,v_{j_k}\}$. Note that $S_G$ and $S_H$
are multisets. Let $\lambda_k(G)=a$, $\lambda_k(H)=b$. Without loss
of generality, $b\geq a\geq\lfloor\frac{k}{2}\rfloor$.

To prove Theorem~\ref{th-Lower}, we use the following five lemmas.
\begin{lemma}\label{lem2-1}
If $|S_G|=|S_H|=k$, we can construct at least $a+b$
edge-disjoint $S$-Steiner trees in $G\Box H$.
\end{lemma}
\begin{proof}
Since $\lambda_k(H)=b$, it follows that there are $b$ edge-disjoint
$S_H$-Steiner trees, say $T_1, T_2,\ldots, T_b$. Without loss of
generality, let $T_1, T_2,\ldots, T_{b'}$ be the edge-disjoint
minimal $S_H$-Steiner trees such that $V(T_i)=S_H$, and $T_{b'+1},
T_{b'+2},$ $\ldots, T_{b}$ be the edge-disjoint minimal
$S_H$-Steiner trees such that $V(T_i)\supset S_H$. From the definition
of $b'$, we have $b'\leq \lfloor k/2\rfloor$. Since $b\geq a\geq
\lfloor k/2\rfloor$, it follows that $a\geq b'$. For each $i \
(1\leq i\leq k)$, let $T_1(u_i), T_2(u_i), \cdots, T_{b}(u_i)$ be
the edge-disjoint $S_H(u_i)$-Steiner trees in $H(u_i)$ corresponding
to $T_1, T_2,\ldots, T_{b}$ in $H$, where
\[
S_H(u_i)=\{(u_i,v_{j_t})\,|\,1\leq t\leq k\}.
\]

Since $\lambda_k(G)=a$, it follows that there are $a$ edge-disjoint
minimal $S_G$-Steiner trees, say $T_1', T_2',\ldots, T_a'$. For each
$j \ (1\leq j\leq k)$, let $T_1'(v_j), T_2'(v_j),$ $ \cdots,
T_{a}'(v_j)$ be the edge-disjoint $S_G(v_j)$-Steiner trees in
$G(v_j)$ corresponding to $T_1', T_2',\ldots, T_{a}'$ in $G$, where
\[
S_G(v_j)=\{(u_{i_s},v_j)\,|\,1\leq s\leq k\}.
\]

\begin{fact}\label{fact1}
For any $S_H(u_{i_s})$-Steiner tree $T_p(u_{i_s}) \ (1\leq s\leq k,
\ 1\leq p\leq b)$ and any $S_G(v_{j_t})$-Steiner tree $T_q'(v_{j_t})
\ (1\leq t\leq k, \ 1\leq q\leq a)$, we can find two edge-disjoint
$S$-Steiner trees in $(\bigcup_{s=1}^kH(u_{i_s}))\cup
(\bigcup_{t=1}^kG(v_{j_t}))$.
\end{fact}
\begin{proof}
Note that $T_q'(v_{j_1})\cup T_p(u_{i_2})\cup T_p(u_{i_3})\cup
\ldots \cup T_p(u_{i_k})$ and $T_p(u_{i_1})\cup T_q'(v_{j_2})\cup
T_q'(v_{j_3})\cup \ldots \cup T_q'(v_{j_k})$ are two edge-disjoint
$S$-Steiner trees in $(\bigcup_{s=1}^kH(u_{i_s}))\cup
(\bigcup_{t=1}^kG(v_{j_t}))$.
\end{proof}

From Fact~\ref{fact1}, we can find $2b'$ edge-disjoint $S$-Steiner
trees in
\[
\left(\bigcup_{p=1}^{b'}\bigcup_{s=1}^kT_p(u_{i_s})\right)\bigcup
\left(\bigcup_{q=1}^{b'}\bigcup_{t=1}^kT_q'(v_{j_t})\right).
\]

Since $b\geq a\geq b'$, from Fact~\ref{fact1}, we can also find
$2a-2b'$ edge-disjoint $S$-Steiner trees in
\[
\left(\bigcup_{p=b'+1}^{a}\bigcup_{s=1}^kT_p(u_{i_s})\right)\bigcup
\left(\bigcup_{q=b'+1}^{a}\bigcup_{t=1}^kT_q'(v_{j_t})\right).
\]
It remains to find $(a+b)-2b'-(2a-2b')=b-a$ edge-disjoint
$S$-Steiner trees in addition to the trees above. Note that we still have
edge-disjoint $S_H(u_{i_s})$-Steiner trees
$T_{a+1}(u_{i_s}),T_{a+2}(u_{i_s}),\ldots,T_{b}(u_{i_s})$ for each
$s \ (1\leq s\leq k)$. Observe that $T_{a+1}$ is a $S_H$-Steiner
tree such that $V(T_{a+1})\supset S_H$. Without loss of generality,
let $V(T_{a+1})-S_H=\{v_{j_{k+1}},v_{j_{k+2}},\ldots,v_{j_{k+h}}\}$.

\begin{figure}[tbp]
\centering
\includegraphics[width=0.80\linewidth]{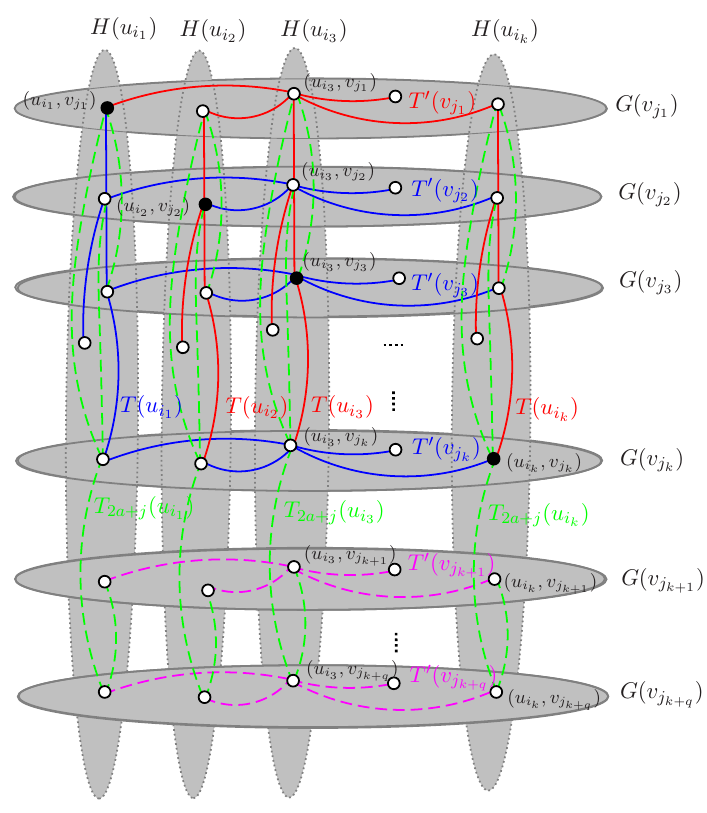}
\caption{Graphs used in the proof of Lemma~\ref{lem2-1}.}
\label{fig:lemma-1}
\end{figure}

Note that there are $a$ edge-disjoint $S_G(v_{j_{k+c}})$-Steiner
trees $T_t'(v_{j_{k+c}})$ in $G(v_{j_{k+c}})$, where $1\leq t\leq a$
and $1\leq c\leq h$. Since $\lambda_k(H)=b$, it follows that there
are $b$ edge-disjoint $S_H$-Steiner trees, say $T_1, T_2,\ldots,
T_b$. Without loss of generality, let $T_1, T_2,\ldots, T_{b'}$ be
the edge-disjoint minimal $S_H$-Steiner trees such that
$V(T_i)=S_H$, and $T_{b'+1}, T_{b'+2},\ldots, T_{b}$ be the
edge-disjoint minimal $S_H$-Steiner trees such that $V(T_i)\supset
S_H$. From the definition of $b'$, we have $b'\leq \lfloor
k/2\rfloor$. Since $b\geq a\geq \lfloor k/2\rfloor$, it follows that
$a\geq b'$. If there is a vertex $v\in V(T_{b'+i}) \ (1\leq i\leq
b-b')$, then the degree of $v$ in $T_{b'+i}$ is at least $2$. Since
there are at most $k$ edges from $v$ to $S_H$, it follows that the
number of edge-disjoint $S_H$-Steiner tree containing $v$ is at most
$\lfloor k/2\rfloor$.

Then we have the following fact.

\begin{fact}\label{fact2}
For any $S_H(u_{i_s})$-Steiner tree $T_p(u_{i_s}) \ (1\leq s\leq k,
\ a+1\leq p\leq b)$ and any $S_G(v_{j_{k+c}})$-Steiner tree
$T_q'(v_{j_{k+c}}) \ (a+1\leq q\leq b, \ 1\leq c\leq h)$, we can
find one $S$-Steiner tree in
\[
\left(\bigcup_{s=1}^{k}H(u_{i_s})\right)\bigcup
\left(\bigcup_{c=1}^{h}G(v_{j_{k+c}})\right).
\]
\end{fact}

From Fact~\ref{fact2}, we can find $(b-a)$ edge-disjoint $S$-Steiner
trees, and the total number of edge-disjoint $S$-Steiner trees is
$a+b$, as desired.
\end{proof}

\begin{lemma}\label{lem2-2}
If $2\leq |S_H|<k$ and $|S_G|=k$, we can construct at
least $a+b$ edge-disjoint $S$-Steiner trees in $G\Box H$.
\end{lemma}
\begin{proof}
Note that $S_G=\{u_{i_1},u_{i_2},\ldots,u_{i_k}\}$, and
$S_H=\{v_{j_1},v_{j_2},\ldots,v_{j_k}\}$. Since $S_H$ is a
multiset, we can assume that $|S_H|=d$ and
$S_H=\{v_{j_1},v_{j_2},\ldots,$ $v_{j_d}\}$. Without loss of
generality, we may let
\[
S\cap G(v_{j_1})=\{(u_{i_c},v_{j_1})\,|\,1\leq c\leq r\}.
\]

For each $s$ with $1\leq s\leq k$, the equality
$\lambda_k(H(u_{i_s}))=b$ yields $b$ edge-disjoint
$S_H(u_{i_s})$-Steiner trees, denoted by
\[
T_1(u_{i_s}),T_2(u_{i_s}),\ldots,T_b(u_{i_s}),
\]
where
\[
S_H(u_{i_s})=\{(u_{i_s},v_{j_t}):1\leq t\leq d\}.
\]
For each $t$ with $1\leq t\leq d$, the equality
$\lambda_k(G(v_{j_t}))=a$ yields $a$ edge-disjoint
$S_G(v_{j_t})$-Steiner trees, denoted by
\[
T'_1(v_{j_t}),T'_2(v_{j_t}),\ldots,T'_a(v_{j_t}),
\]
where
\[
S_G(v_{j_t})=\{(u_{i_s},v_{j_t}):1\leq s\leq k\}.
\]

\begin{fact}\label{fact3}
For any $S_H(u_{i_s})$-Steiner tree $T_p(u_{i_s}) \ (1\leq s\leq k,
\ 1\leq p\leq b)$ and any $S_G(v_{j_t})$-Steiner tree $T_q'(v_{j_t})
\ (1\leq t\leq d, \ 1\leq q\leq a)$, we can find two edge-disjoint
$S$-Steiner trees in $(\bigcup_{s=1}^kH(u_{i_s}))\cup
(\bigcup_{t=1}^dG(v_{j_t}))$.
\end{fact}
\begin{proof}
Note that $T_q'(v_{j_1})\cup T_p(u_{i_{r+1}})\cup
T_p(u_{i_{r+2}})\cup \cdots \cup T_p(u_{i_k})$ and
$T_q'(v_{j_2})\cup T_q'(v_{j_3})\cup \cdots \cup T_q'(v_{j_d})\cup
T_p(u_{i_1})\cup T_p(u_{i_2})\cup \cdots \cup T_p(u_{i_r})$ are two
edge-disjoint Steiner trees in 8
$\left(\bigcup_{s=1}^kH(u_{i_s})\right)\cup
\left(\bigcup_{t=1}^dG(v_{j_t})\right)$.
\end{proof}

From Fact~\ref{fact3}, we can find $2a$ edge-disjoint $S$-Steiner
trees in
\[
\left(\bigcup_{p=1}^{a}\bigcup_{s=1}^kT_p(u_{i_s})\right)\cup
\left(\bigcup_{q=1}^{b}\bigcup_{t=1}^kT_q'(v_{j_t})\right).
\]
It remains to find $(a+b)-2a=b-a$ edge-disjoint $S$-Steiner trees
in addition to the trees above. Note that we still have edge-disjoint
$S_H(u_{i_s})$-Steiner trees
$T_{a+1}(u_{i_s}),T_{a+2}(u_{i_s}),\ldots,T_{b}(u_{i_s})$ for each
$s \ (1\leq s\leq k)$. Observe that $T_{a+1}$ is a $S_H$-Steiner
tree such that $V(T_{a+1})\supset S_H$. Without loss of generality,
let $V(T_{a+1})-S_H=\{v_{j_{k+1}},v_{j_{k+2}},\ldots,v_{j_{k+c}}\}$.

Note that there are $a$ edge-disjoint $S_G(v_{j_{k+e}})$-Steiner
trees $T_q'(v_{j_{k+e}})$ in $G(v_{j_{k+e}})$, where $1\leq q\leq a$
and $1\leq e\leq c$. From Fact~\ref{fact2}, we can find $(b-a)$
edge-disjoint $S$-Steiner trees, and the total number of
edge-disjoint $S$-Steiner trees is $a+b$, as desired.
\end{proof}

\begin{lemma}\label{lem2-3}
If $2\leq |S_G|<k$ and $|S_H|=k$, we can construct at
least $a+b$ edge-disjoint $S$-Steiner trees in $G\Box H$.
\end{lemma}

\begin{proof}
Note that $S_H=\{v_{j_1},v_{j_2},\ldots,v_{j_k}\}$, and
$S_G=\{u_{i_1},u_{i_2},\ldots,u_{i_k}\}$. Since $S_G$ is a
multiset, we can assume that $|S_G|=d$ and
$S_G=\{u_{i_1},u_{i_2},\ldots,u_{i_d}\}$. Without loss of
generality, we may let
\[
S\cap H(u_{i_1})=\{(u_{i_1},v_{j_c})\,|\,1\leq c\leq r\}.
\]

For each $s$ with $1\leq s\leq d$, the equality
$\lambda_k(H(u_{i_s}))=b$ yields $b$ edge-disjoint
$S_H(u_{i_s})$-Steiner trees, denoted by
\[
T_1(u_{i_s}),T_2(u_{i_s}),\ldots,T_b(u_{i_s}),
\]
where
\[
S_H(u_{i_s})=\{(u_{i_s},v_{j_t}):1\leq t\leq k\}.
\]
For each $t$ with $1\leq t\leq k$, the equality
$\lambda_k(G(v_{j_t}))=a$ yields $a$ edge-disjoint
$S_G(v_{j_t})$-Steiner trees, denoted by
\[
T'_1(v_{j_t}),T'_2(v_{j_t}),\ldots,T'_a(v_{j_t}),
\]
where
\[
S_G(v_{j_t})=\{(u_{i_s},v_{j_t}):1\leq s\leq d\}.
\]

\begin{fact}\label{fact4}
For any $S_H(u_{i_s})$-Steiner tree $T_p(u_{i_s}) \ (1\leq s\leq d,
\ 1\leq p\leq b)$ and any $S_G(v_{j_t})$-Steiner tree $T_q'(v_{j_t})
\ (1\leq t\leq k, \ 1\leq q\leq a)$, we can find two edge-disjoint
$S$-Steiner trees in $(\bigcup_{s=1}^dH(u_{i_s}))\cup
(\bigcup_{t=1}^kG(v_{j_t}))$.
\end{fact}
\begin{proof}
The subgraphs
\[
T_p(u_{i_1})\cup T_q'(v_{j_{r+1}})\cup T_q'(v_{j_{r+2}})
 \cup\cdots\cup T_q'(v_{j_k})
\]
and
\[
T_p(u_{i_2})\cup T_p(u_{i_3})\cup\cdots\cup T_p(u_{i_d})
 \cup T_q'(v_{j_1})\cup T_q'(v_{j_2})\cup\cdots\cup T_q'(v_{j_r})
\]
are two edge-disjoint $S$-Steiner trees in
\[
\left(\bigcup_{s=1}^d H(u_{i_s})\right)
 \cup \left(\bigcup_{t=1}^k G(v_{j_t})\right).
\]
\end{proof}

From Fact~\ref{fact4}, we can find $2a$ edge-disjoint $S$-Steiner
trees in
\[
\left(\bigcup_{p=1}^{a}\bigcup_{s=1}^kT_p(u_{i_s})\right)\cup
\left(\bigcup_{q=1}^{b}\bigcup_{t=1}^kT_q'(v_{j_t})\right).
\]
It remains to find $(a+b)-2a=b-a$ edge-disjoint $S$-Steiner trees
in addition to the trees above. Note that we still have edge-disjoint
$S_H(u_{i_s})$-Steiner trees
$T_{a+1}(u_{i_s}),T_{a+2}(u_{i_s}),\ldots,T_{b}(u_{i_s})$ for each
$s \ (1\leq s\leq k)$. Observe that $T_{a+1}$ is a $S_H$-Steiner
tree such that $V(T_{a+1})\supset S_H$. Without loss of generality,
let $V(T_{a+1})-S_H=\{v_{j_{k+1}},v_{j_{k+2}},\ldots,v_{j_{k+c}}\}$.

Note that there are $a$ edge-disjoint $S_G(v_{j_{k+e}})$-Steiner
trees $T_q'(v_{j_{k+e}})$ in $G(v_{j_{k+e}})$, where $1\leq q\leq a$
and $1\leq e\leq c$. From Fact~\ref{fact2}, we can find $(b-a)$
edge-disjoint $S$-Steiner trees, and the total number of
edge-disjoint $S$-Steiner trees is $a+b$, as desired.
\end{proof}

\begin{figure}[tbp]
\centering
\includegraphics[width=4in]{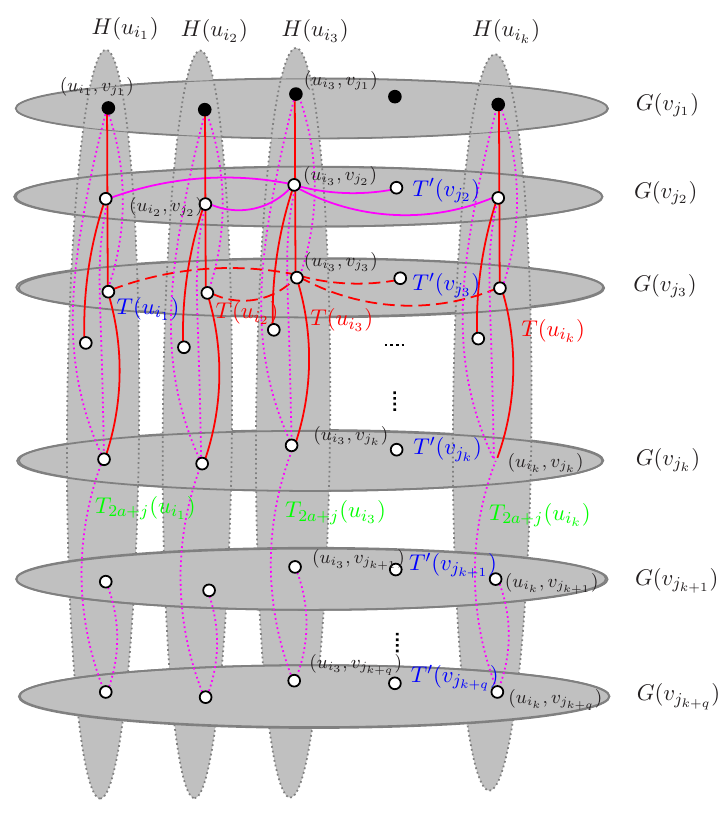}
\caption{Graphs used in the proof of Lemma~\ref{lem2-4}.}\label{fig:lemma-4}
\end{figure}

\begin{lemma}\label{lem2-4}
If $|S_G|=k$ and $|S_H|=1$ or $|S_G|=1$ and $|S_H|=k$,
we can construct $a+b$ edge-disjoint $S$-Steiner trees in $G\Box H$.
\end{lemma}
\begin{proof}
Without loss of generality, let $|S_G|=k$ and $|S_H|=1$.  Thus,
\[
S_G=\{u_{i_1},u_{i_2},\ldots,u_{i_k}\}
\quad\text{and}\quad
S_H=\{v_{j_1}\},
\]
where the latter equality concerns the set of distinct vertices in the
multiset $\{v_{j_1},v_{j_2},\ldots,v_{j_k}\}$. Since
$\lambda_k(G(v_{j_1}))=a$, there are $a$ edge-disjoint
$S_G(v_{j_1})$-Steiner trees, and they are also $a$ edge-disjoint
$S$-Steiner trees.

Let $S_H^*=\{v_{j_1},v_{r_2},\ldots,v_{r_k}\}$, where
$v_{r_2},\ldots,v_{r_k}\in V(H)-v_{j_1}$. Since $\lambda_k(H)=b$, it
follows that there are $b$ edge-disjoint $S_H^*$-Steiner trees, say
$T_1, T_2,\ldots, T_b$. Then there are $b$ edge-disjoint $S$-Steiner
trees in
\[
\bigcup_{p=1}^{b}\left[\left(\bigcup_{q=1}^{a}\bigcup_{v_i\in
V(T_p)-v_{j_1}}T_q'(v_i)\right)\bigcup
\left(\bigcup_{i=1}^{k}T_p(u_i)\right)\right],
\]
where $T_q'(v_i)$ and $T_p(u_i)$ are as defined in Lemma~\ref{lem2-1}.
Thus, we can construct $a+b$ edge-disjoint $S$-Steiner trees in $G\Box
H$, as desired.
\end{proof}

\begin{lemma}\label{lem2-5}
If $|S_G|<k$ and $|S_H|<k$, we can construct $a+b$
edge-disjoint $S$-Steiner trees.
\end{lemma}
\begin{proof}
Let $|S_G|=c$ and $|S_H|=d$. Note that
$S_G=\{u_{i_1},u_{i_2},\ldots,u_{i_k}\}$ and
$S_H=\{v_{j_1},v_{j_2},\ldots,v_{j_k}\}$ are both multisets. We
assume that $S_G=\{u_{i_1},u_{i_2},\ldots,u_{i_c}\}$ and
$S_H=\{v_{j_1},v_{j_2},\ldots,v_{j_d}\}$. Let
\begin{align*}
S_G^*={}&\{u_{i_1},u_{i_2},\ldots,u_{i_c},
          u_{i_{c+1}}^*,u_{i_{c+2}}^*,\ldots,u_{i_k}^*\},\\
S_H^*={}&\{v_{j_1},v_{j_2},\ldots,v_{j_d},
          v_{j_{d+1}}^*,v_{j_{d+2}}^*,\ldots,v_{j_k}^*\},
\end{align*}
where
\[
u_{i_{c+1}}^*,u_{i_{c+2}}^*,\ldots,u_{i_k}^*
 \in V(G)\setminus\{u_{i_1},u_{i_2},\ldots,u_{i_c}\}
\]
and
\[
v_{j_{d+1}}^*,v_{j_{d+2}}^*,\ldots,v_{j_k}^*
 \in V(H)\setminus\{v_{j_1},v_{j_2},\ldots,v_{j_d}\}.
\]

For each $t$ with $1\leq t\leq d$, the equality
$\lambda_k(G(v_{j_t}))=a$ yields $a$ edge-disjoint
$S_G(v_{j_t})$-Steiner trees in $G(v_{j_t})$, denoted by
\[
T'_1(v_{j_t}),T'_2(v_{j_t}),\ldots,T'_a(v_{j_t}),
\]
where
\[
S_G(v_{j_t})=
\{(u_{i_s},v_{j_t}):1\leq s\leq c\}
\cup\{(u_{i_s}^*,v_{j_t}):c+1\leq s\leq k\}.
\]
For each $s$ with $1\leq s\leq c$, the equality
$\lambda_k(H(u_{i_s}))=b$ yields $b$ edge-disjoint
$S_H(u_{i_s})$-Steiner trees, denoted by
\[
T_1(u_{i_s}),T_2(u_{i_s}),\ldots,T_b(u_{i_s}),
\]
where
\[
S_H(u_{i_s})=
\{(u_{i_s},v_{j_t}):1\leq t\leq d\}
\cup\{(u_{i_s},v_{j_t}^*):d+1\leq t\leq k\}.
\]

\begin{fact}\label{fact5}
$(1)$ For any $S_H(u_{i_s})$-Steiner tree $T_p(u_{i_s}^*) \ (c+1\leq
s\leq k)$ and any $S_G(v_{j_t})$-Steiner tree $T_q'(v_{j_t}) \
(1\leq t\leq c, \ 1\leq q\leq a)$, we can find one $S$-Steiner tree
in $T_p(u_{i_s}^*)\cup (\bigcup_{t=1}^cT_q'(v_{j_t}))$ for any $q$.

$(2)$ For any $S_G(v_{j_t})$-Steiner tree $T_q'(v_{j_t}^*) \
(d+1\leq t\leq k)$ and any $S_H(u_{i_s})$-Steiner tree $T_p(u_{i_s})
\ (1\leq s\leq d, \ 1\leq p\leq b)$, we can find one $S$-Steiner
tree in
$\bigcup_{v_j\in V(T_p)\setminus V(S_H)}T_q'(v_j)
 \cup(\bigcup_{s=1}^dT_p(u_{i_s}))$, for any $p$.
\end{fact}

\begin{figure}[tbp]
\centering
\includegraphics[width=4in]{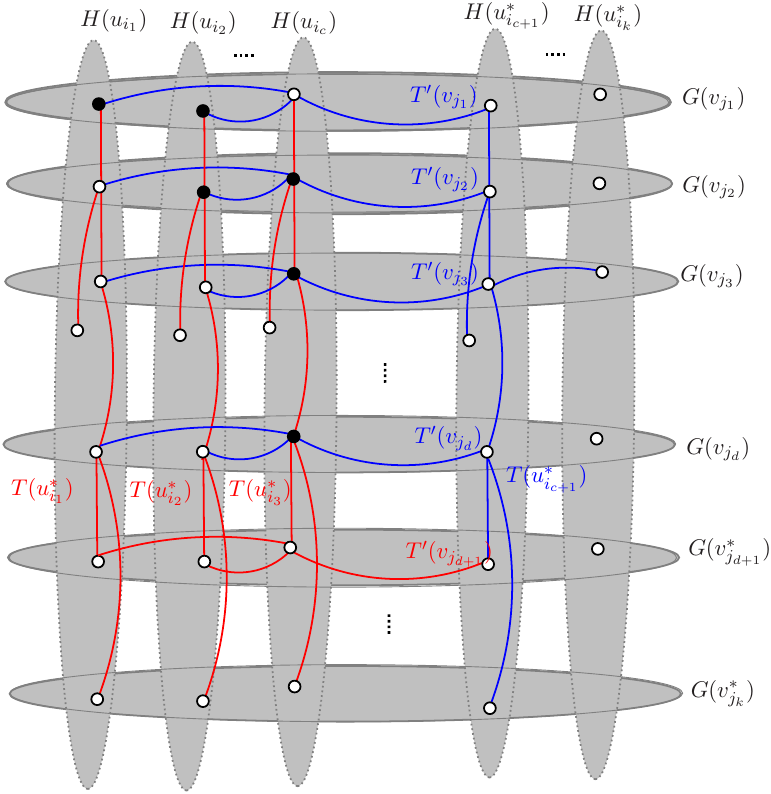}
\caption{Graphs used in the proof of Lemma~\ref{lem2-5}.}\label{fig:lemma-5}
\end{figure}

From $(1)$ of Fact~\ref{fact5}, we can find $a$ edge-disjoint
$S$-Steiner trees. From $(2)$ of Fact~\ref{fact5}, we can find $b$
edge-disjoint $S$-Steiner trees. Note that all the $S$-Steiner trees
are edge-disjoint, as desired.
\end{proof}

By Lemmas~\ref{lem2-1}--\ref{lem2-5}, we have \[\lambda_k(G\Box H)\geq
\lambda_k(G)+\lambda_k(H).\]

To show the sharpness of Proposition~\ref{pro-Upper} and
Theorem~\ref{th-Lower}, we consider the following example.

\begin{example}
Let $\mathcal{F}$ be a graph class containing graphs $F_i$ obtained
from a complete graph $K_n$ and a vertex $u$ by adding $i$ edges
between $u$ and $K_n$, where $\lfloor\frac{k}{2}\rfloor\leq i\leq
n-1$ and $k\leq n-1-\lceil k/2\rceil$. Choose $G,H\in \mathcal{F}$.
Then $\lambda_k(G)=\delta(G)=i$ and $\lambda_k(H)=\delta(H)=i$. From
Proposition~\ref{pro-Upper} and Theorem~\ref{th-Lower}, we have
\[
\lambda_k(G)+\lambda_k(H)\leq \lambda_k(G\Box H)\leq
\delta(G)+\delta(H),
\]
and hence $\lambda_k(G\Box H)=\lambda_k(G)+\lambda_k(H)$.
\end{example}

\section{Results for networks}

A \emph{two-dimensional grid graph} is an $m\times n$ graph
Cartesian product $P_n\,\Box\,P_m$ of path graphs on $m$ and $n$
vertices. For more details on grid graphs, we refer to \cite{Calkin,
ItaiRodeh}.

\begin{proposition}
Let $m$ and $n$ be integers with $m\geq n\geq 3$. Then
\[
\lambda_k(P_m\Box P_n)= \begin{cases}
2 & \text{ if $3\leq k<\min\{n,m\}$,}\\
1 & \text{ if $\lceil \frac{2mn-m-n+2}{2}\rceil <k\leq mn$,}\\
1\text{ or }2 & \text{ if $\min\{n,m\}\leq k\leq \lceil
\frac{2mn-m-n+2}{2}\rceil$.}
\end{cases}
\]
\end{proposition}
\begin{proof}
Suppose that $k<\min\{n,m\}$. For $1\leq j\leq n$, let
\[
e(P_j)=\{(u_i,v_j)(u_{i+1},v_j)\,|\,1\leq i\leq m-1\},
\]
and for $1\leq i\leq m$, let
\[
e(Q_i)=\{(u_i,v_j)(u_{i},v_{j+1})\,|\,1\leq j\leq n-1\}.
\]
For any $S\subseteq V(P_m\Box P_n)$ with $|S|=k$, since $k<m$, there is an index $j$ such that
$S\cap V(P_j)=\varnothing$, where $V(P_j)=\{(u_i,v_j)\,|\,1\leq i\leq m\}$.
Since $k<n$, there is an index $i$ such that
$S\cap V(Q_i)=\varnothing$, where $V(Q_i)=\{(u_i,v_j)\,|\,1\leq
j\leq n\}$. The subgraph induced by the edge set
\[
\left(\bigcup_{a=1}^{j-1}e(P_a)\right)\cup
\left(\bigcup_{a=j+1}^{n}e(P_a)\right)\cup e(Q_i)
\]
contains an $S$-Steiner tree. Similarly, the subgraph induced by
\[
\left(\bigcup_{b=1}^{i-1}e(Q_b)\right)\cup
\left(\bigcup_{b=i+1}^{m}e(Q_b)\right)\cup e(P_j)
\]
contains an $S$-Steiner tree. Since the two $S$-Steiner trees are
edge-disjoint, it follows that
$\lambda_k(P_m\Box P_n)\geq 2$. Since $\lambda_k(P_m\Box P_n)\leq
\delta(P_m\Box P_n)=2$, it follows that $\lambda_k(P_m\Box P_n)= 2$.

Suppose that $\lceil \frac{2mn-m-n+2}{2}\rceil <k\leq mn$. It is
clear that $\lambda_k(P_m\Box P_n)\geq 1$. We will show that
$\lambda_k(P_m\Box P_n)=1$. For any $S\subseteq V(P_m\Box P_n)$ with
$|S|=k$, if we want to find two edge-disjoint $S$-Steiner trees,
then we need at least $2k-2$ edges. Since $k>\lceil
\frac{2mn-m-n+2}{2}\rceil$, it follows that $2k-2>2mn-m-n=e(P_m\Box
P_n)$, a contradiction.
\end{proof}

An \emph{torus} is the Cartesian product of two cycles $C_{m},C_{n}$
of size at least three. The two cycles need not have the same size.

\begin{proposition}
Let $m$ and $n$ be integers with $m\geq n\geq 3$. Then
\[
\lambda_k(C_m\Box C_n)= \begin{cases}
2\text{ or }3 & \text{ if $3\leq k<\lceil \frac{2mn+3}{3}\rceil$,}\\
2 & \text{ if $\lceil \frac{2mn+3}{3}\rceil <k\leq mn$.}
\end{cases}
\]
\end{proposition}
\begin{proof}
Suppose that $\lceil \frac{2mn+3}{3}\rceil <k\leq mn$. It is clear
that $\lambda_k(C_m\Box C_n)\geq \lambda_{mn}(C_m\Box C_n)\geq 2$.
We will show that $\lambda_k(C_m\Box C_n)=2$. For any $S\subseteq
V(C_m\Box C_n)$ with $|S|=k$, if we want to find three edge-disjoint
$S$-Steiner trees, then we need at least $3k-3$ edges. Since
$k>\lceil \frac{2mn+3}{3}\rceil$, it follows that
$3k-3>2mn=e(C_m\Box C_n)$, a contradiction.

Suppose that $3\leq k<\lceil \frac{2mn+3}{3}\rceil$. Clearly,
$\lambda_k(C_m\Box C_n)\geq \lambda_{mn}(C_m\Box C_n)\geq 2$. Since
there are two adjacent vertices of degree $4$, we have
$\lambda_k(C_m\Box C_n)\leq \delta(C_m\Box C_n)-1=3$.
\end{proof}

\section{Concluding remarks}

We give a lower bound of $\lambda_k(G\Box H)$ under the condition
$\lambda_k(H)\geq \lambda_k(G)\geq \lfloor \frac{k}{2}\rfloor$. The
case that $\lambda_k(G)\leq \lfloor \frac{k}{2}\rfloor$ is still
open. It is also open to determine $\lambda_k(P_m\Box P_n)=1$ or $2$
for $\min\{n,m\}\leq k\leq \lceil \frac{2mn-m-n+2}{2}\rceil$;
$\lambda_k(C_m\Box C_n)=2$ or $3$ for $3\leq k<\lceil
\frac{2mn+3}{3}\rceil$.

\acknowledgments

We would like to thank the anonymous referees for a number of
helpful comments and suggestions. This work was partially supported
by the National Natural Science Foundation of China (Nos. 11601254,
11551001, 11871280, 12271259), the National Key R\&D Program of
China (No. 2020YFB2104004), Qinghai Key R\&D and Transformation
Projects (No. 2021-GX-112), and the Qinghai Key Laboratory of
Internet of Things Project (No. 2017-ZJ-Y21).

\end{document}